\documentclass[11pt,reqno]{amsart}
\usepackage{amsmath,amssymb, amsrefs}
\usepackage[margin=2.5cm]{geometry}
\usepackage{graphicx}
\usepackage{hyperref}

\theoremstyle{plain}
\newtheorem{thm}{Theorem}[section]

\newtheorem{que}[thm]{Question}

\newtheorem*{mresult}{Main Result}

\theoremstyle{definition}

\newtheorem*{rem}{Remark}

\theoremstyle{definition}

\numberwithin{equation}{section}

\newcommand{\vir}{\mathrm{vir}}

\newcommand{\ov}{\overline}
\newcommand{\lb}{\left[}
\newcommand{\rb}{\right]}
\newcommand{\pp}{\mathbb{P}}
\newcommand{\zz}{\mathbb{Z}}
\renewcommand{\qq}{\mathbb{Q}}
\renewcommand{\L}{\xi_X}
\newcommand{\R}{\mathcal{R}}
\newcommand{\pr}{\mathrm{pr}}

\newcommand{\ev}{\mathrm{ev}}

\usepackage{hyperref}
\hypersetup{
	colorlinks,
	citecolor=red,
	filecolor=black,
	linkcolor=blue,
	urlcolor=black
}

\def \hf{\hspace*{0.5cm}}
\begin{document}
\title{Elliptic Gromov-Witten Invariants of Del-Pezzo surfaces}
\author[C. Chaudhuri]{Chitrabhanu Chaudhuri}
\address{School of Mathematics, IISER Pune}
\email{chitrabhanu@iiserpune.ac.in}
\author[N. Das]{Nilkantha Das}
\address{School of Mathematics, National Institute of Science 
         Education and Research, Bhubaneswar (HBNI), Odisha 752050, India}
\email{nilkantha.das@niser.ac.in } 

\begin{abstract}
We obtain a formula for the number of genus one curves with a variable  
complex structure of a given degree on a del-Pezzo surface that pass 
through an appropriate number of generic points of the surface. 
This is done using Getzler's relationship among cohomology classes
of certain codimension 2 cycles in $\ov{M}_{1,4}$ 
and 
recursively computing the genus one Gromov-Witten 
invariants of del-Pezzo surfaces. 
Using completely different methods, this problem has been solved earlier by 
Bertram and Abramovich (\cite{AB}), 
Ravi Vakil (\cite{Vak}), Dubrovin and Zhang (\cite{Dub}) and more  recently 
using Tropical geometric methods by 
M.~Shoval and E.~Shustin (\cite{Sh}). We also subject our formula to several low degree 
checks and compare them to the numbers obtained by the earlier authors. 

\end{abstract}

\maketitle
\tableofcontents

\section{Introduction}
One of the most fundamental problems in enumerative algebraic geometry is: 
\begin{que}
What is $E^{(g)}_d$, the number of genus $g$ degree $d$ curves in $\mathbb{CP}^2$ (with a variable complex structure) that pass through $3d-1+g$ generic points?  
\end{que}
Although the computation of $E^{(g)}_d$ 
is a classical question, a complete solution to the above problem (even for genus zero) 
was unknown until the early
$90^{' \textnormal{s}}$ when Ruan--Tian (\cite{RT}) and Kontsevich--Manin (\cite{KoMa})
obtained a formula for $E^{(0)}_d$. \\ 
\hf \hf The computation of $E^{(g)}_d$ is now very well understood from several different perspectives. 
The formula by Caporasso--Harris \cite{CH}, computes $E^{(g)}_d$ for all $g$ and $d$. Since then, the computation of 
$E^{(g)}_d$ has been studied from many different perspectives; these include (among others) the algorithm 
by Gathman (\cite{AnGa1}, \cite{AnGa2}) and the  
method of virtual localization by Graber and Pandharipande (\cite{GP}) to compute the genus $g$ Gromov-Witten invariants of $\mathbb{CP}^n$ 
(although for $n>2$ and $g>0$, the Gromov-Witten invariants are not enumerative). 
More recently, the problem of computing $E^{(g)}_d$ has been studied using the method of tropical geometry by Mikhalkin in \cite{Mi} (using the results 
of that paper, one can in principle compute $E^{(g)}_d$ for all $g$ and $d$). \\ 
\hf \hf A more  general situation is as follows: let $X$ be a projective manifold and 
and $\beta \in H_2(X;\mathbb{Z})$ a given homology class. 
Given cohomology classes $\mu_1, \ldots, \mu_k \in H^*(X, \qq)$, the $k$-pointed genus $g$ 
Gromov-Witten invariant of $X$ is defined to be
\begin{equation} \label{gw}
  N^{(g)}_{\beta,X}(\mu_1, \ldots, \mu_k) :=  
  \int \limits_{\ov{M}_{g, k}(X, \beta)} 
  \ev_1^*(\mu_1) \smile \ldots \smile \ev_{k}^*(\mu_k) \smile 
  ~\lb \ov{M}_{g, k}(X, \beta) \rb ^{\vir}, 
\end{equation}
where $\ov{M}_{g,k}(X, \beta)$ denotes the moduli space of genus $g$ stable maps
into $X$ with $k$ marked points representing $\beta$ and $\textnormal{ev}_i$ denotes 
the $i^{\textnormal{th}}$ evaluation map. For $g=0$, this is a smooth, irreducible 
and proper Deligne-Mumford stack and has a fundamental class. However, for $g>0$,
$\ov{M}_{g,k}(X,\beta)$ is not smooth or irreducible, hence it does not posses 
a fundamental class. Behrend, Behrend-Fantechi and Li-Tian, have however defined 
the virtual fundamental class  
\begin{equation*}
  \lb \ov{M}_{g, k}(X, \beta) \rb ^{\vir} \in H^{2 \Theta}(\ov{M}_{g,k}(X,\beta)), 
  \qquad \Theta := c_1(TX)\cdot \beta + (\dim X -3)(1-g) + k;
\end{equation*}
which is used to define the Gromov-Witten invariants (see \cite{B},\cite{BF} and
\cite{LiTi}).
When all the $\mu_1, \ldots, \mu_k$ represent the class Poincare dual to a point 
(and the degree of the cohomology class that is being paired in \eqref{gw}, is equal 
to the virtual dimension of the moduli space), then we abbreviate 
$N^{(g)}_{\beta, X}(\mu_1,\ldots, \mu_k)$ as $N^{(g)}_{\beta}$. 
The number of genus g curves of degree $\beta$ in $X$, that pass through $c_1(TX) \cdot \beta + (\dim X -3)(1-g)$ generic points is denoted by $E^{(g)}_{\beta}$. 
In general, $E^{(g)}_{\beta}$ is 
not necessarily equal to $N^{(g)}_{\beta}$, i.e. the Gromov-Witten invariant is not 
necessarily enumerative (this happens for example when $X:= \mathbb{CP}^3$ and $g=1$). \\ 
\hf \hf An important class of surfaces for which the enumerative geometry 
is particularly 
important 
are Fano surfaces, which are also called del-Pezzo surfaces 
(see section \ref{dp_defn} for the definition of a del-Pezzo surface). 
When $g=0$, it is proved in (\cite{Pandh_Gott}, Theorem 4.1, Lemma 4.10) 
that for del-Pezzo surfaces $N^{(0)}_{\beta} = E^{(0)}_{\beta}$.\\ 
\hf \hf In \cite{Vak}, Vakil generalizes the approach of Caporasso-Harris in \cite{CH} to compute the numbers $E^{(g)}_{\beta}$ for 
all $g$ and $\beta$ for del-Pezzo surfaces. It is also shown in (\cite{Vak}, Section 4.2) that all the genus $g$ Gromov-Witten invariants of del-Pezzo surfaces 
are enumerative (i.e. $N^{(g)}_{\beta} = E^{(g)}_{\beta}$). The enumerative geometry of del-Pezzo surfaces has also been studied extensively  
by Abramovich and Bertram (in \cite{AB}). 
More recently, this question has been approached using methods of tropical geometry. 
In \cite{Sh}, M.~Shoval and E.~Shustin  
give a formula to compute all the genus $g$ Gromov-Witten invariants of del-Pezzo surfaces using methods of tropical geometry.\\ 
\hf \hf The genus one Gromov-Witten invariants of $\mathbb{CP}^n$ can also be computed from a completely different method 
from the ones developed in \cite{AnGa1}, \cite{AnGa2} and \cite{GP}. In \cite{EG}, Getzler finds a relationship among certain codimension 
two cycles in $\overline{M}_{1,4}$ and uses that to compute the genus one Gromov-Witten invariants of $\mathbb{CP}^2$ and $\mathbb{CP}^3$. 
In \cite{EGH}, using ideas from Physics, Eguchi, Hori and Xiong made a remarkable conjecture concerning the genus $g$ Gromov-Witten invariants of 
projective manifolds; this is known as the Virasoro conjecture. The conjecture in particular produces an explicit formula for $N^{(1)}_d$ 
(for $\mathbb{CP}^2$), which aprori looks very different from the formula obtained by Getzler (in \cite{EG}).   
It is shown by Pandharipande (in \cite{Pan}), 
that the formula obtained by Getzler for $\mathbb{CP}^2$ is equivalent to a completely different looking formula 
predicted in \cite{EGH}. \\ 
\hf \hf In this paper,  we extend the approach of Getzler to compute the genus one Gromov-Witten invariants of del-Pezzo surfaces. 
The formula we obtain has a completely different appearance from the one obtained by Vakil in \cite{Vak}. 
We verify that 
our final numbers are consistent with the numbers he obtains (see section \ref{ld_check} for details).   \\ 
\hf \hf The Virasoro conjecture for projective manifolds (which is conjectured in \cite{EGH}) has been a topic of active research in 
mathematics for the last twenty years. In \cite{Dub}, 
Dubrovin and Zhang  
compute the genus one Gromov-Witten invariants of $\mathbb{CP}^1 \times \mathbb{CP}^1$ by showing that it follows 
from the Virasoro conjecture. We have verified that our numbers agree with all the numbers computed by them (\cite{Dub}, Page 463). 
They prove that the genus zero and genus one Virasoro Conjecture is true for all projective manifolds having semi-simple quantum cohomology. 
It is proved in \cite{BaMa} that the quantum cohomology of del-Pezzo surfaces is semi simple. It would be interesting to see if one 
can use the result of this paper and 
apply the method of  
\cite{Pan} to  
obtain a formula for the genus one Gromov-Witten invariants of del-Pezzo surfaces, analogous to 
the one predicted for $\mathbb{CP}^2$ by Eguchi, Hori and Xiong (in \cite{EGH}). That would give a \textit{direct} confirmation of the 
Virasoro conjecture in genus one for del-Pezzo surfaces. 
A detailed survey of the Virasosro conjecture is given in \cite{EGV}. 

\section{Main Result}
The main result of this paper is the following:
\begin{mresult} \label{main_result}
Let $X$ be a del-Pezzo surface and $\beta\in H_2(X, \mathbb{Z})$ be a given effective
homology class. We obtain a formula for $N^{(1)}_{\beta}$ (equation \eqref{rec_formula}) 
using Getzler's relation 
\end{mresult}

\begin{rem}
We note that by (\cite{Vak}, Section 4.2), we conclude that  $N^{(1)}_{\beta} = E^{(1)}_{\beta}$. 
Alternatively, we note that $N^{(1)}_{\beta} = E^{(1)}_{\beta}$ follows from (\cite{Zi_Red}, Theorem 1.1).  
\end{rem}

Our formula for $N^{(1)}_{\beta}$ is a recursive formula, involving $N^{(0)}_{\beta}$. 
The latter can be computed via the algorithm given in \cite{KoMa} and \cite{Pandh_Gott}. 
The base case of our recursive formula are given by equations \eqref{rec_blowup} and 
\eqref{rec_p1xp1}. We have written a $\mbox{C++}$ program 
that implements \eqref{rec_formula}; it is available on our web page: 
\begin{center}
\url{http://www.iiserpune.ac.in/~chitrabhanu/}.
\end{center}


\section{Recursive formula} 
We will now give the recursive formula to compute $N^{(1)}_{\beta}$.
First, we will develop some notation that is used throughout this paper. Let
\begin{equation*}
  \begin{array}{lll}
    \L & := c_1(T X), \qquad & \text{for both the cohomology class and the divisor}, \\
    \kappa_{\beta} & := \L \cdot \beta, \qquad  & \text{where } \beta \in H_2(X,\zz), \\
    b_2(X) & := \dim H_2(X, \qq), \qquad & \text{the second betti number of } X,\\
    d_X & := \L \cdot \L, \qquad & \text{the degree of } X.
  \end{array}  
\end{equation*}
Moreover, $\cdot$ is used for both the cup product in cohomology as well as cap product
between a homology and a cohomology class.

We are now ready to state the formula. First, let us define the following four quantities:  
\begin{alignat*}{3}
  T_1 := & \sum_{\beta_1 + \beta_2 + \beta_3 = \beta} &
         & \binom{\kappa_{\beta}-2}{\kappa_{\beta_2} - 1, \kappa_{\beta_3} - 1} 
           2 \kappa_{\beta_2} \kappa_{\beta_3}^2 (\beta_1 \cdot \beta_2) \\ 
         & & & \bigg( \big( 4 \kappa_{\beta_1} + \kappa_{\beta_2} 
         - 2\kappa_{\beta_3} \big) (\beta_2 \cdot \beta_3) 
         - 3\kappa_{\beta_2}(\beta_1 \cdot \beta_3) \bigg) 
         N^{(1)}_{\beta_1} N^{(0)}_{\beta_2} N^{(0)}_{\beta_3}, \\ 
  \\
  T_2 := & \sum_{\beta_1 + \beta_2 = \beta}  \bigg[ &
         & \binom{\kappa_{\beta}-2}{\kappa_{\beta_1}-1} 
           4 \kappa_{\beta_2}^2 \bigg( 2\kappa_{\beta_1}\kappa_{\beta_2} -                     
           \kappa_{\beta_2}^2 - 3 d_X (\beta_1 \cdot \beta_2) \bigg) + \\
         & & & \binom{\kappa_{\beta}-2}{\kappa_{\beta_1}} 
           2\kappa_{\beta_2}\bigg( d_X (\beta_1 \cdot \beta_2)
           \big(4\kappa_{\beta_1} + \kappa_{\beta_2} \big) +
           2 \kappa_{\beta_1}\kappa_{\beta_2} \big( 2 \kappa_{\beta_1}
           - \kappa_{\beta_2} \big) \bigg) \ \bigg] 
           N^{(1)}_{\beta_1} N^{(0)}_{\beta_2}, \\ 
\end{alignat*}
\begin{alignat*}{3}
  T_3 := & - \frac{1}{12} \sum_{\beta_1 + \beta_2 = \beta} 
           \binom{\kappa_{\beta}-2}{\kappa_{\beta_1}-1} 
           \kappa_{\beta_2}^2 (\beta_1 \cdot \beta_2)\bigg[ &
         &\kappa_{\beta_1}^2\bigg( \kappa_{\beta_1}
           - 2\kappa_{\beta_2}  - 6(\beta_1 \cdot \beta_2) \bigg) \\            
  &      & & + \kappa_{\beta_2}(\beta_1 \cdot \beta_1)
           \bigg( 4\kappa_{\beta_1} + \kappa_{\beta_2}) \bigg) \bigg] 
         N^{(0)}_{\beta_1} N^{(0)}_{\beta_2}, \\ 
  T_4 := & - \frac{1}{12} \kappa_{\beta}^3 
           \bigg( (2+ b_2(X))\kappa_{\beta}  
           - d_X \bigg) N_{\beta}^{(0)}.
\end{alignat*}
The number $N^{(1)}_{\beta}$ satisfies the following recursive relation: 
\begin{equation} \label{rec_formula}
  6 d_X^2 N^{(1)}_{\beta} = T_1 + T_2 + T_3 + T_4. 
\end{equation}
We will  now give the initial conditions for the recursion \eqref{rec_formula}.   
Let $X$ be $\mathbb{P}^2$ blown up at upto $k=8$ points. 
Then the initial condition of the recursion is 
\begin{equation} \label{rec_blowup}
  N^{(1)}_L = 0 \qquad  \textnormal{and} \qquad 
  N^{(1)}_{E_i} = 0 \qquad \forall i = 1 ~~\textnormal{to} ~~k. 
\end{equation}
Here $L$ denotes the class of a line and $E_i$ denotes the exceptional divisors.  
If $X:= \mathbb{P}^1 \times \mathbb{P}^1$, then 
\begin{equation} \label{rec_p1xp1}
  N^{(1)}_{e_1} = 0 \qquad \textnormal{and} \qquad N^{(1)}_{e_2} =0. 
\end{equation}
Here $e_1$ and $e_2$ denote the class of $[\textnormal{pt} \times \mathbb{P}^1]$ 
and $[\mathbb{P}^1 \times \textnormal{pt}]$ 
respectively. The initial conditions \eqref{rec_blowup} and \eqref{rec_p1xp1}, 
combined with the values of $N^{(0)}_{\beta}$ obtained from \cite{KoMa} and 
\cite{Pandh_Gott}, give us the values of $N^{(1)}_{\beta}$ for any $\beta$. 

\begin{rem}
  We would like to mention that the formula \eqref{rec_formula} yields Getzler's
  recursion relation, equation $(0.1)$ of \cite{EG}, after some symmetrization of 
  the summation indices of $T_1$ and $T_3$.
\end{rem}  

\section{Del-Pezzo surfaces}
\label{dp_defn}
A del-Pezzo surface $X$ is a smooth projective algebraic surface with an ample 
anti-canonical divisor $\L$. The degree of the surface is defined to be 
the self-intersection number  
\begin{equation*}
  d_X = \L \cdot \L.
\end{equation*}
This degree $d_X$ varies between $1$ and $9$. $X$ can be obtained as a blow-up 
of $\pp^2$ at $k=9-d_X$ general points, except, when $d_X = 8$ the surface can also 
be $\pp^1 \times \pp^1$. 

If $X$ has degree $9-k$ and is not $\pp^1 \times \pp^1$, then we have the blow up
morphism $Bl: X \to \pp^2$. We denote by $E_1, \ldots, E_k$ the exceptional divisors 
of $Bl$ and by $L$ the pull-back of the class of a hyperplane in $\pp^2$. We have
\begin{equation*}
  H^2(X, \zz) = \zz \langle L, E_1, \ldots, E_k \rangle,
\end{equation*}
and $L \cdot L = 1$, $E_i \cdot E_i = -1$, $L \cdot E_i = E_i \cdot E_j = 0$ 
for all $i,j \in \{1,\ldots,k\}$ with $i \neq j$. The anti-canonical divisor is 
given by $\L = 3L - E_1 -\ldots -E_k$.

If $X = \pp^1 \times \pp^1$, let $e_1 = \pr_1^* [\textnormal{pt}]$ and $e_2 = \pr_2^*[\textnormal{pt}]$, 
then $\L = 2 e_1 + 2 e_2$, $e_1 \cdot e_2 = 1$ and $e_2 \cdot e_1 =1$ whereas $e_i \cdot e_i = 0$
for $i = 1,2$.

\section{Basic Strategy}

\begin{center}
\begin{figure}[h]
\begin{tabular}{ccc}
  \includegraphics[scale=0.35]{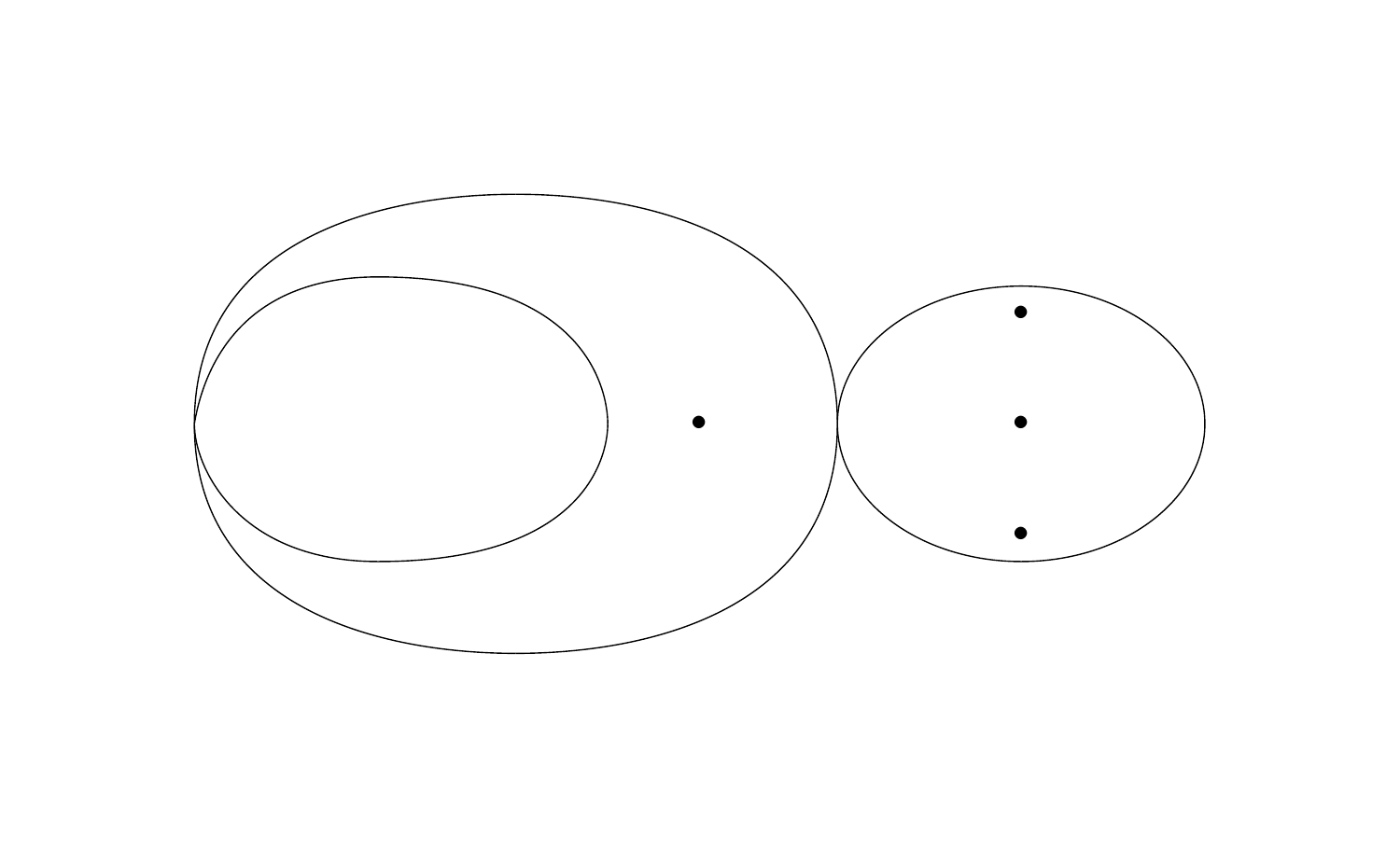}
  & 
  \includegraphics[scale=0.35]{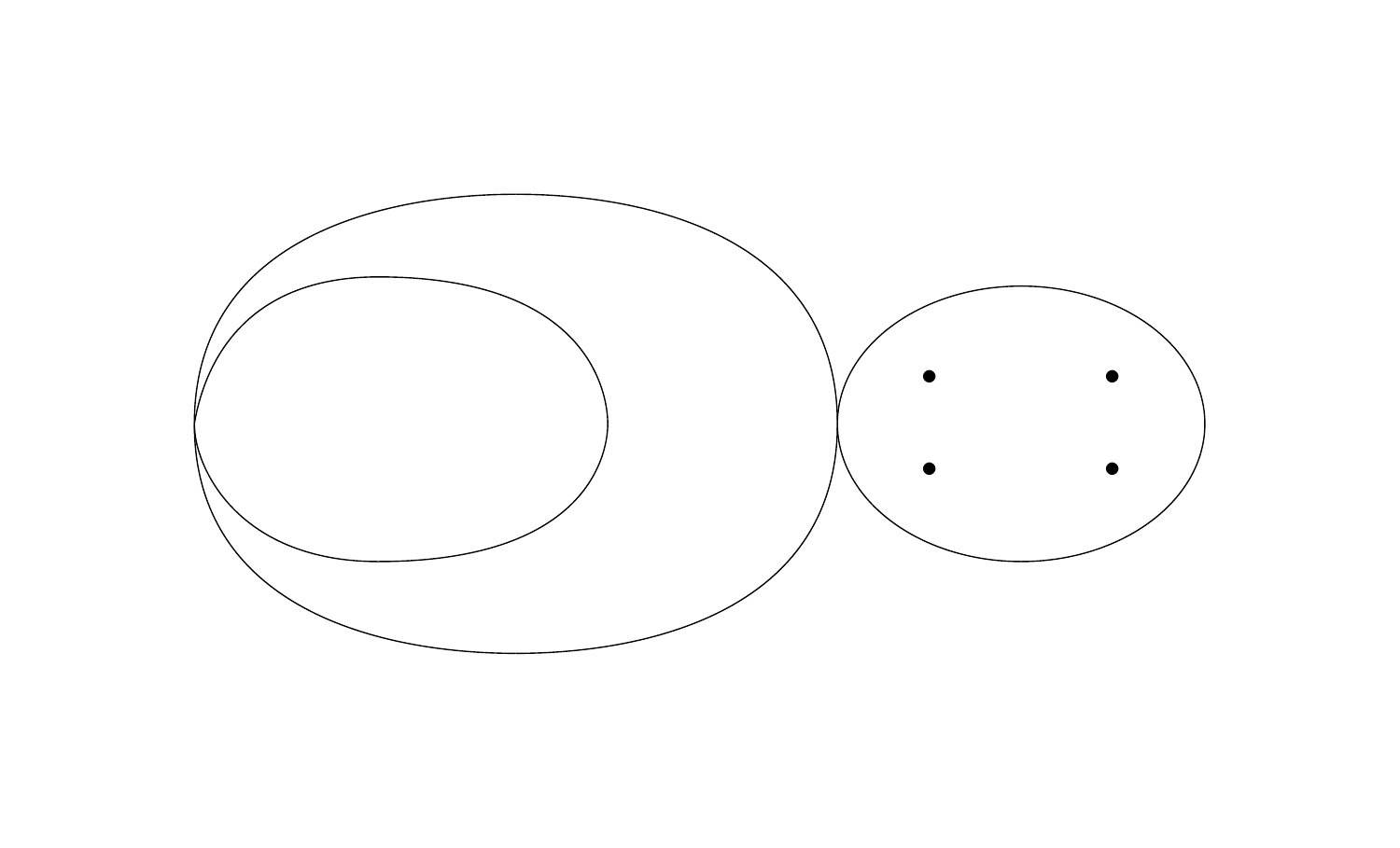}
  &
  \includegraphics[scale=0.35]{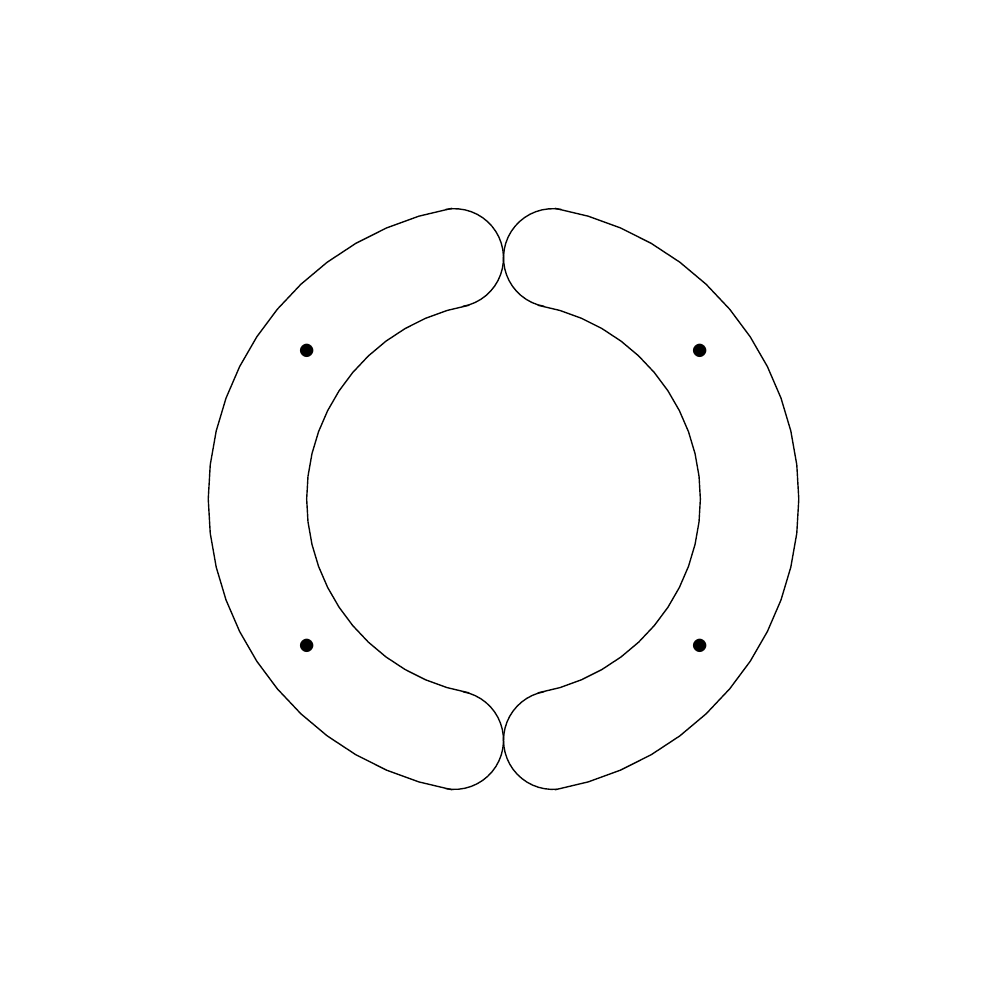}
  \\
  $\delta_{0,3}$ & $\delta_{0,4}$ & $\delta_{0,0}$ \\  
  \includegraphics[scale=0.3]{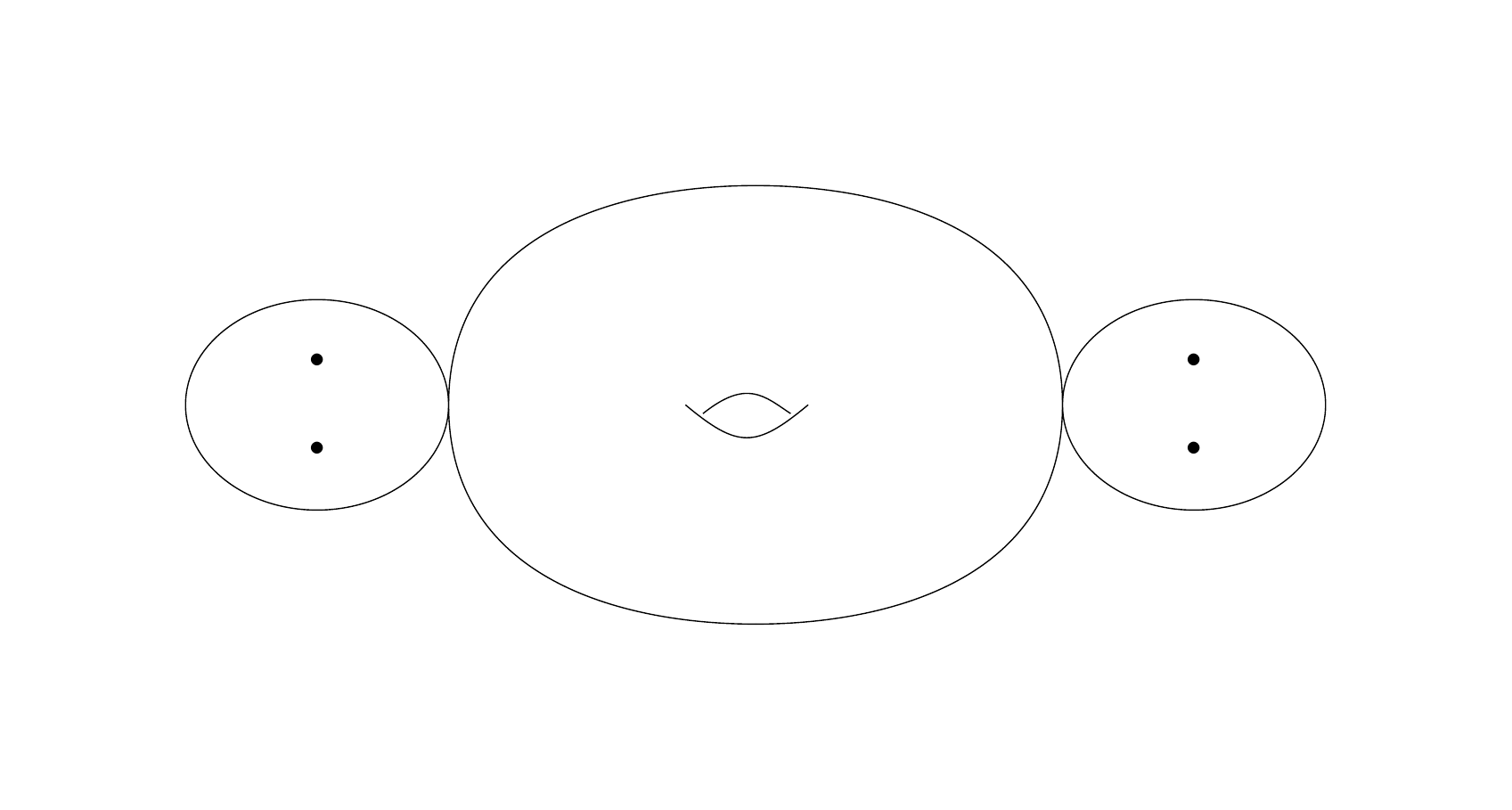}
  & 
  \includegraphics[scale=0.3]{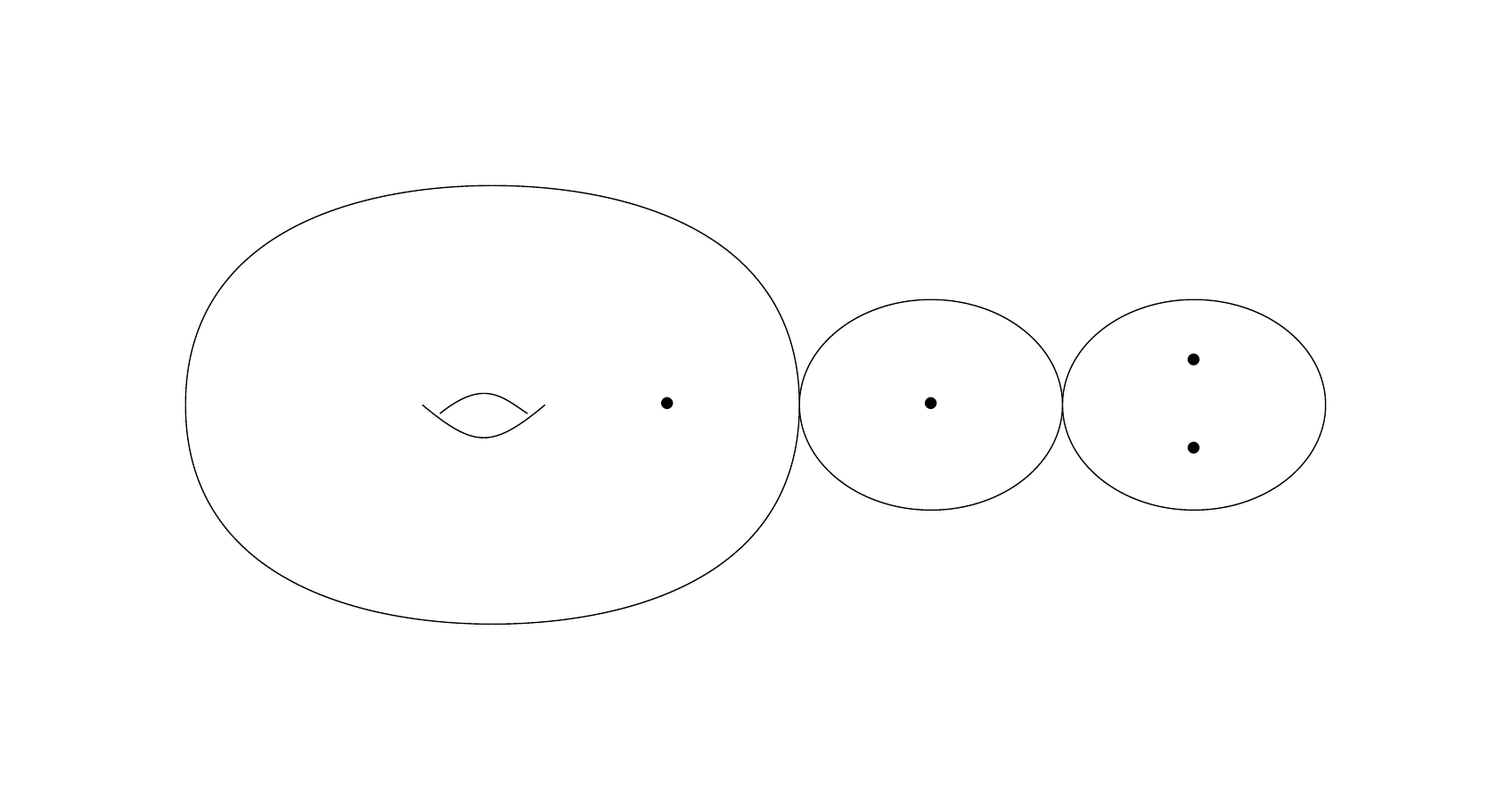}
  &
  \includegraphics[scale=0.3]{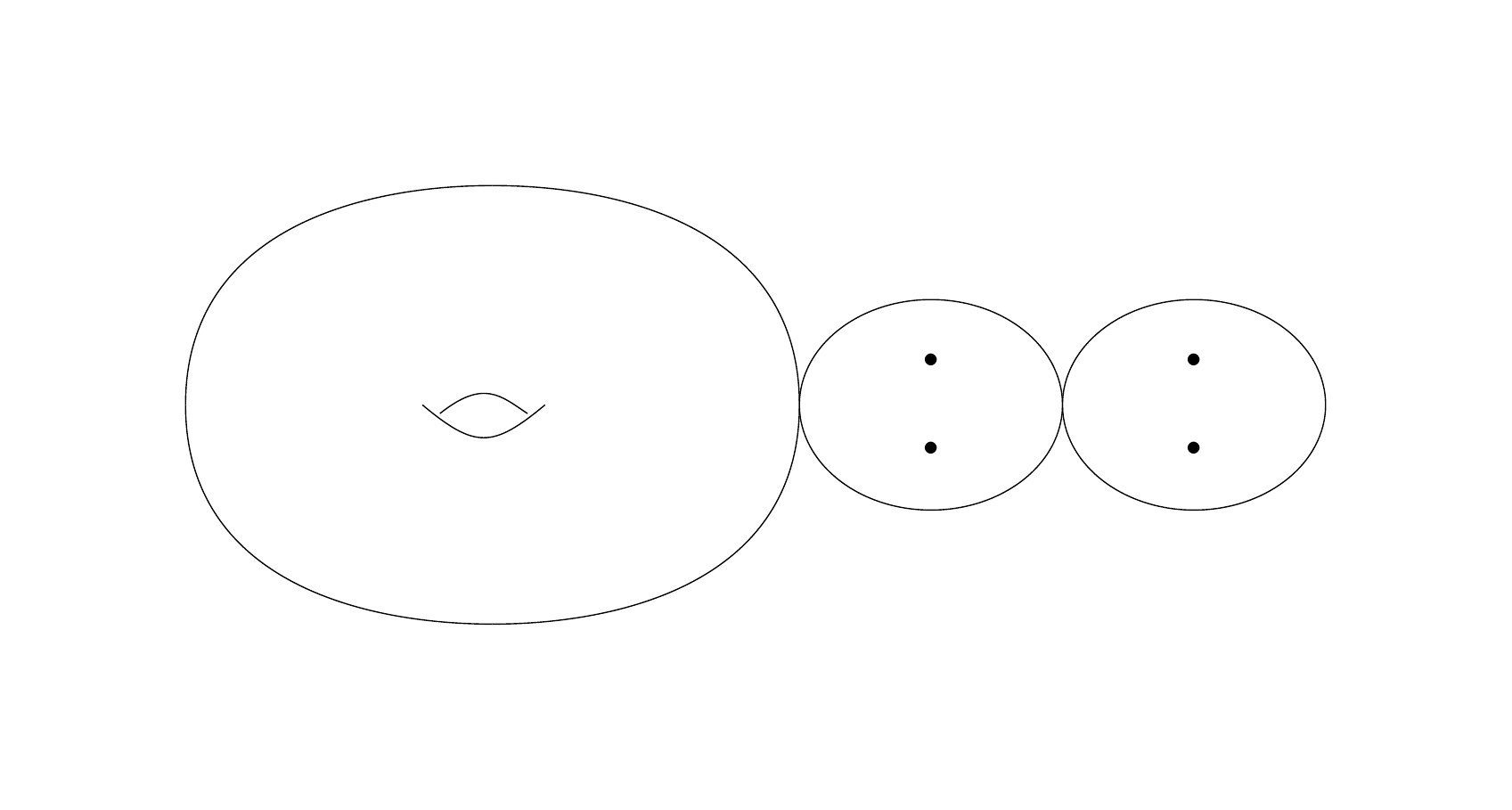}\\
  $\delta_{2,2}$ & $\delta_{2,3}$ & $\delta_{2,4}$ \\
  &  \includegraphics[scale=0.3]{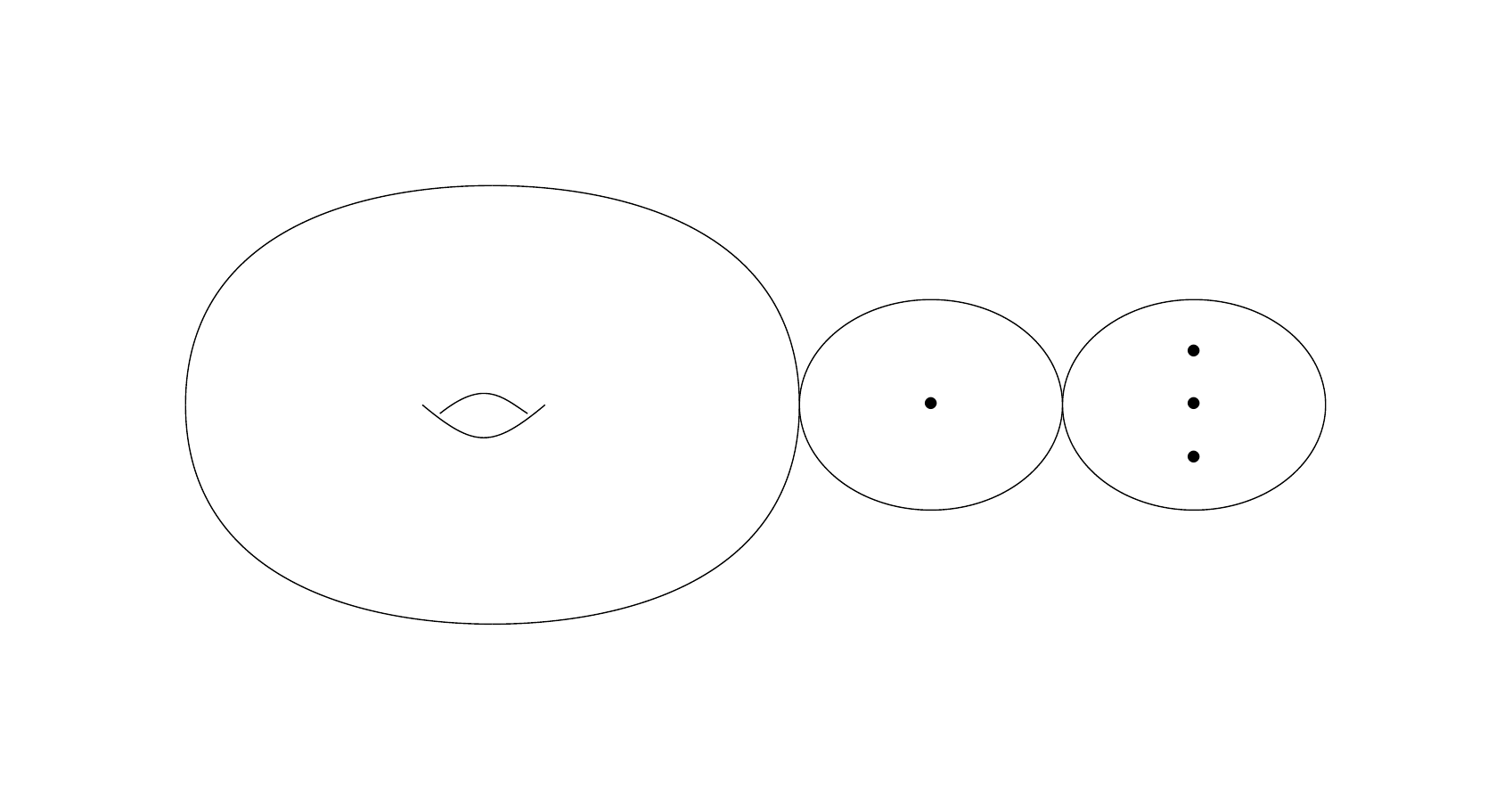} &  \\
  & $\delta_{3,4}$ & \\
\end{tabular} 
\caption{Codimension 2 strata in $\ov{M}_{1,4}$.}
 \label{fig:strata}
\end{figure}
\end{center}

We will now recall the basic setup of \cite{EG}, where Getzler computes the number 
$N^{(1)}_d$ when $X$ is $\mathbb{CP}^2$. First, let us consider the space $\ov{M}_{1,4}$, 
the moduli space of genus one curves with four marked  points. We shall be interested 
in certain $S_4$ invariant codimension 2 boundary strata in $\ov{M}_{1,4}$ which we 
list in Figure \ref{fig:strata}. In the figure we draw the topological type and the 
marked point distribution of the generic curve in each strata. We use the same 
nomenclature as \cite{EG} except for $\delta_{0,0}$ which was denoted by $\delta_{\beta}$
in \cite{EG}, (to avoid confusion between notations). See section 1 of \cite{EG} for a 
list of all the codimension 2 strata. There the strata are denoted by the dual graph 
of the generic curve.

These strata define cycles in $H^4(\ov{M}_{1,4},\qq)$. Let us now define the following 
cycle in $H^4(\overline{M}_{1,4}, \mathbb{Q})$, given by 
\begin{equation*}
  \R  := -2 \delta_{2,2}+\frac{2}{3} \delta_{2,3}  
         + \frac{1}{3} \delta_{2,4} - \delta_{3,4} 
         -\frac{1}{6} \delta_{0,3} -\frac{1}{6} \delta_{0,4} 
         + \frac{1}{3} \delta_{0,0}.  
\end{equation*}
The main result of \cite{EG} is that $\R =0$. This will subsequently be referred to as 
Getzler's relation. In \cite{Pan}, Pandharipande has shown that this relation, in fact, 
comes from a rational equivalence.

Now we explain how to obtain our formula. Consider the natural forgetful morphism
\begin{equation*}
  \pi: \overline{M}_{1,\kappa_{\beta} + 2}(X, \beta) \longrightarrow \overline{M}_{1,4}.
\end{equation*}
We shall pull-back the cycle $\R$ to $H^*(\ov{M}_{1,\kappa_{\beta} + 2}(X, \beta), \qq)$ 
and intersect it with a cycle of complementary dimension; that will give us an equality 
of numbers and subsequently the formula. Let $\mu \in H^4(X,\qq)$ be the class of a point. 
Define 
\begin{equation*}
  \mathcal{Z} := \ev_1^*(\L) \cdot \ldots \cdot \ev_4^*(\L) \cdot 
                 \ev_5^*(\mu) \cdot \ldots \cdot \ev_{\kappa_{\beta}+2}^*(\mu).
\end{equation*}
The class $\L$ is used since it is ample and hence numerically effective.
Since $\R =0$ by Getzler's relation, we conclude that
\begin{equation} \label{eq:relation}
  \int\limits_{\ov{M}_{1,\kappa_{\beta} + 2}(X, \beta)} 
      (\pi^*\R \cdot \mathcal{Z}) 
      \cdot \lb \ov{M}_{1,\kappa_{\beta} + 2}(X, \beta)\rb^{\vir} 
      = 0. 
\end{equation}
We can also compute the left hand side of \eqref{eq:relation} using the composition axiom
for Gromov-Witten invariants which will give us the recursive formula. 

\section{Axioms for Gromov-Witten Invariants} 
We shall make use of certain axioms for Gromov-Witten invariants. These are quite standard,
see for example \cite{CoxKatz}, however for completeness we list them here. We assume 
$X$ is a smooth projective variety.

\begin{description}
\item[Degree axiom] If $\deg \mu_1 + \ldots + \deg \mu_n \neq 2n + 
     2\kappa_{\beta} + 2(3- \dim X)(g - 1)$ then
     \begin{equation*}
       N^{(g)}_{\beta, X}(\mu_1,\ldots,\mu_n) = 0. 
     \end{equation*}
\item[Fundamental class axiom] If $[X]$ is the fundamental class of $X$
     and $2g + n \geq 4$ or $\beta \neq 0$, then 
     \begin{equation*}
       N^{(g)}_{\beta, X}([X], \mu_1,\ldots, \mu_{n-1}) = 0.
     \end{equation*}     
\item[Divisor axiom] If $D$ is a divisor of $X$ and $2g + n \geq 4$. then 
     \begin{equation*}
       N^{(g)}_{\beta, X}(D, \mu_1,\ldots, \mu_{n-1}) = 
       (D \cdot \beta) N^{(g)}_{\beta, X}(\mu_1,\ldots, \mu_{n-1}).
     \end{equation*}         
\item[Composition axiom] This is a bit complicated to write down, so we refer
     to \cite{EG}, section 2.11. It is a combination of the splitting and 
     reduction axioms of \cite{KoMa} section 2.
\end{description}

We also need the following results which do not follow from the above axioms:
\begin{equation*} \label{eq:gen0}
  N^{(0)}_{0,X}(\mu_1, \mu_2, \mu_3) = \int_X \mu_1 \smile \mu_2 \smile \mu_3,
\end{equation*}
and
\begin{equation*} \label{eq:gen1}
  N^{(1)}_{0,X}(\mu) = - \frac{1}{24} c_1(TX) \cdot \mu.
\end{equation*}

\section{Intersection of cycles} 
Now we are in a position to compute the left hand side of \eqref{eq:relation}. Fix 
a homogeneous basis $\{\gamma_1, \ldots, \gamma_{b(X)}\}$ of $H^*(X,\qq)$. Let
$g_{ij} = \int_X \gamma_i \smile \gamma_j$ and $((g^{ij})) = ((g_{ij}))^{-1}$. 
For a cycle $\delta$ in $H^*(\ov{M}_{g,n}(X,\beta),\qq)$, we introduce the following 
notation
\begin{equation*}
  N^{\delta}_{\beta,X}(\mu_1,\ldots,\mu_n) = 
  \int\limits_{\ov{M}_{g,n}(X, \beta)}
  \delta \cdot \ev_1^*(\mu_1) \cdots \ev_n^*(\mu_n) \cdot 
  \lb \ov{M}_{g,n}(X, \beta)\rb^{\vir}.
\end{equation*}

Let $\mu_1 = \ldots = \mu_4 = \xi_X$, and $\mu_5 = \ldots = 
\mu_{\kappa_{\beta}+2} = [pt]$ be the class of a point. If $\delta = \pi^*  
\delta_{2,2}$, by the composition axiom
\begin{alignat*}{3}
  N^{\delta}_{\beta,X} =
    & N^{\delta}_{\beta,X}(\mu_1,\ldots,\mu_{\kappa_{\beta}+2}) && \\ 
  = & \sum_{\substack{\beta_1 + \beta_2 + \beta_3 = \beta \\ A, B, C}} \ 
    \sum_{i,j,k,l} g^{ij} g^{kl}
    && N^{(1)}_{\beta_1,X}(\gamma_i, \gamma_k, \mu_{\alpha} | \alpha \in A)\\ 
  & && \times N^{(0)}_{\beta_2,X}(\gamma_j, \mu_{\alpha} | \alpha \in B) 
      \times N^{(0)}_{\beta_3,X}(\gamma_l, \mu_{\alpha} | \alpha \in C),
\end{alignat*}
where the second sum is over $i,j,k,l$ ranging from $1$ to $b(X)$ and
the first sum is over disjoint sets $A, B, C$ satisfying
\begin{equation*}
  A \sqcup B \sqcup C = \{1,\ldots, \kappa_{\beta}+2\}, \quad
  |B \cap \{1,2,3,4\}| = |C \cap \{1,2,3,4\}| =2.
\end{equation*}
Note that if $\beta_1 , \beta_2, \beta_3 > 0$, by the degree axiom the only 
non-trivial terms occur when $|A| = \kappa_{\beta_1}, |B| = \kappa_{\beta_2}+1, 
|C| = \kappa_{\beta_3}+1$. The limiting case $\beta_1 =0$ does not yield anything, 
however $\beta_2 = 0$ or $\beta_3 = 0$ have non-trivial contributions to the sum. 
When $\beta_3 =0, \beta_1, \beta_2 > 0$, the non-trivial contribution occurs precisely 
when $|C| = 2, \gamma_l = [X]$, $|A| = \kappa_{\beta_1}-1, \gamma_k = [pt]$, and 
$|B| = \kappa_{\beta_2} +1$. Finally when $\beta_2 = \beta_3 = 0$, the only non-zero
term occurs when $|B|=|C|=2$, $\gamma_l = \gamma_j = [X]$ and $\gamma_k = \gamma_i
= [pt]$. Making use of the fact that for any $\sigma, \tau \in H^*(X,\qq)$
\begin{equation*}
  \sum_{i = 1}^{b(X)} \sum_{j=1}^{b(X)} g^{ij} (\sigma \cdot \gamma_i)
  (\gamma_j \cdot \tau) = (\sigma \cdot \tau),
\end{equation*}
we obtain the following expression
\begin{align}
  N^{\pi^*\delta_{2,2}}_{\beta,X}  =
    & 3 (\L \cdot \L)^2 N^{(1)}_{\beta}  \nonumber \\
  & + 3 \sum_{\beta_1 + \beta_2 + \beta_3 = \beta} 
    \binom{\kappa_{\beta}-2}{\kappa_{\beta_2}-1, \kappa_{\beta_3}-1} 
    (\beta_2 \cdot \L)^2 (\beta_3 \cdot \L)^2 
    (\beta_1 \cdot \beta_2)(\beta_1 \cdot \beta_3) 
    N^{(1)}_{\beta_1} N^{(0)}_{\beta_2} N^{(0)}_{\beta_3} \nonumber  \\
  & + 6 \sum_{ \beta_1 + \beta_2 = \beta} 
    \binom{\kappa_{\beta}-2}{\kappa_{\beta_1}-1} 
    (\L \cdot \L) (\beta_1 \cdot \beta_2) (\beta_2 \cdot \L)^2 
    N^{(1)}_{\beta_1} N^{(0)}_{\beta_2}. 
\end{align}

Next, let us consider the cycle $\delta_{2,3}$. We then have
\begin{align*}
  N^{\pi^* \delta_{2,3}}_{\beta,X} =
    \sum_{\substack{\beta_1 + \beta_2 + \beta_3 = \beta \\ A, B, C}} \ 
    \sum_{i,j,k,l} g^{ij} g^{kl}
    & N^{(1)}_{\beta_1,X}(\gamma_i, \mu_{\alpha} | \alpha \in A)\\ 
  & \times N^{(0)}_{\beta_2,X}(\gamma_j,\gamma_k, \mu_{\alpha} | \alpha \in B) 
      \times N^{(0)}_{\beta_3,X}(\gamma_l, \mu_{\alpha} | \alpha \in C),
\end{align*}
where the sum is over sets $A, B, C$ satisfying
\begin{equation*}
  A \sqcup B \sqcup C = \{1,\ldots, \kappa_{\beta}+2\}, \quad
  |A \cap \{1,2,3,4\}| = |B \cap \{1,2,3,4\}| = 1.
\end{equation*}
All the cases are similar to the previous calculation except, when $\beta_2 = 0$.
In this case we can either have $|B| = 1, |A| = \kappa_{\beta_1}$, $\gamma_i = [pt]$
and $\gamma_j = [X]$; or $|B| = 1, |C| = \kappa_{\beta_3}$, $\gamma_k = [X]$
and $\gamma_l = [pt]$. We get
\begin{align}
  N^{\pi^* \delta_{2,3}}_{\beta,X} = 
    & 12 \sum_{ \beta_1 + \beta_2 + \beta_3 = \beta}
    \binom{\kappa_{\beta}-2}{\kappa_{\beta_2}-1, \kappa_{\beta_3}-1} 
    (\beta_1 \cdot \L)(\beta_2 \cdot \L) (\beta_3 \cdot \L)^2 
    (\beta_1 \cdot \beta_2)(\beta_2 \cdot \beta_3) 
    N^{(1)}_{\beta_1} N^{(0)}_{\beta_2} N^{(0)}_{\beta_3} \nonumber \\
  & + 12 \sum_{ \beta_1 + \beta_2= \beta} 
    \binom{\kappa_{\beta}-2}{\kappa_{\beta_1}} 
    (\beta_1 \cdot \L) (\beta_2 \cdot \L) 
    \bigg( (\L \cdot \L) (\beta_1 \cdot \beta_2) + 
           (\beta_1 \cdot \L)(\beta_2 \cdot \L) \bigg) 
    N^{(1)}_{\beta_1} N^{(0)}_{\beta_2} \nonumber \\
  & +12 \sum\limits_{ \beta_1 + \beta_2= \beta} 
    \binom{\kappa_{\beta}-2}{\kappa_{\beta_1}-1} 
    (\beta_1 \cdot \L) (\beta_2 \cdot \L)^3  
    N^{(1)}_{\beta_1} N^{(0)}_{\beta_2}. 
\end{align}

Moving on to $\delta_{2,4}$ we have
\begin{align*}
  N^{\pi^* \delta_{2,4}}_{\beta,X} =
    \sum_{\substack{\beta_1 + \beta_2 + \beta_3 = \beta \\ A, B, C}} \ 
    \sum_{i,j,k,l} g^{ij} g^{kl}
    & N^{(1)}_{\beta_1,X}(\gamma_i, \mu_{\alpha} | \alpha \in A)\\ 
  & \times N^{(0)}_{\beta_2,X}(\gamma_j,\gamma_k, \mu_{\alpha} | \alpha \in B) 
      \times N^{(0)}_{\beta_3,X}(\gamma_l, \mu_{\alpha} | \alpha \in C),
\end{align*}
where the sum is over sets $A, B, C$ satisfying
\begin{equation*}
  A \sqcup B \sqcup C = \{1,\ldots, \kappa_{\beta}+2\}, \quad
  |B \cap \{1,2,3,4\}| = |C \cap \{1,2,3,4\}| = 2.
\end{equation*}
Now there is no contribution when $\beta_2 = 0$, however we have a non-trivial 
contribution when $\beta_1 = 0$. We can use \eqref{eq:gen1} to calculate this
\begin{align}
  N^{\pi^* \delta_{2,4}}_{\beta,X} = & \ 
    6 \sum_{ \beta_1 + \beta_2 + \beta_3 = \beta}
    \binom{\kappa_{\beta}-2}{\kappa_{\beta_2}-1, \kappa_{\beta_3}-1} 
    (\beta_2 \cdot \L)^2 (\beta_3 \cdot \L)^2 (\beta_1 \cdot \beta_2)
    (\beta_2 \cdot \beta_3) N^{(1)}_{\beta_1} 
    N^{(0)}_{\beta_2} N^{(0)}_{\beta_3} \nonumber \\   
  & + 6 \sum_{ \beta_1 + \beta_2= \beta} 
    \binom{\kappa_{\beta}-2}{\kappa_{\beta_1}} 
    (\beta_2 \cdot \L)^2 (\L \cdot \L) (\beta_1 \cdot \beta_2) 
    N^{(1)}_{\beta_1} N^{(0)}_{\beta_2} \nonumber \\
  & + 6 \sum_{ \beta_1 + \beta_2= \beta} 
    \left(-\frac{1}{24}\right) \binom{\kappa_{\beta}-2}{\kappa_{\beta_1} -1} 
    (\L \cdot \beta_1)^3 (\beta_2 \cdot \L)^2 (\beta_1 \cdot \beta_2) 
    N^{(0)}_{\beta_1} N^{(0)}_{\beta_2} \nonumber \\
  & + 6 \left(-\frac{1}{24}\right) (\L \cdot \beta)^3 (\L \cdot \L) 
    N^{(0)}_{\beta}.
\end{align}

For $\delta_{3,4}$ we have
\begin{align*}
  N^{\pi^* \delta_{3,4}}_{\beta,X} =
    \sum_{\substack{\beta_1 + \beta_2 + \beta_3 = \beta \\ A, B, C}} \ 
    \sum_{i,j,k,l} g^{ij} g^{kl}
    & N^{(1)}_{\beta_1,X}(\gamma_i, \mu_{\alpha} | \alpha \in A)\\ 
  & \times N^{(0)}_{\beta_2,X}(\gamma_j,\gamma_k, \mu_{\alpha} | \alpha \in B) 
      \times N^{(0)}_{\beta_3,X}(\gamma_l, \mu_{\alpha} | \alpha \in C),
\end{align*}
where the first sum is over sets $A, B, C$ satisfying
\begin{equation*}
  A \sqcup B \sqcup C = \{1,\ldots, \kappa_{\beta}+2\}, \quad
  |B \cap \{1,2,3,4\}| =1, |C \cap \{1,2,3,4\}| = 3.
\end{equation*}
The calculation is similar to the previous cases, so we omit the details. We obtain
\begin{align}
  N^{\pi^* \delta_{3,4}}_{\beta,X} = & \ 
    4  \sum_{ \beta_1 + \beta_2 + \beta_3 = \beta}
    \binom{\kappa_{\beta}-2}{\kappa_{\beta_2} - 1, \kappa_{\beta_3} - 1} 
    (\beta_2 \cdot \L) (\beta_3 \cdot \L)^3 
    (\beta_1 \cdot \beta_2)(\beta_2 \cdot \beta_3) 
    N^{(1)}_{\beta_1} N^{(0)}_{\beta_2} N^{(0)}_{\beta_3} \nonumber \\
  & +4 \sum_{ \beta_1 + \beta_2= \beta} 
    \binom{\kappa_{\beta}-2}{\kappa_{\beta_1}} 
    (\beta_2 \cdot \L)^3 (\beta_1 \cdot \L) 
    N^{(1)}_{\beta_1} N^{(0)}_{\beta_2} \nonumber \\
  & +4 \sum_{ \beta_1 + \beta_2= \beta} 
    \binom{\kappa_{\beta}- 2}{\kappa_{\beta_1} - 1} (\beta_2 \cdot \L)^4 
    N^{(1)}_{\beta_1} N^{(0)}_{\beta_2} \nonumber \\
  & +4 \sum_{ \beta_1 + \beta_2= \beta} 
    \binom{\kappa_{\beta}-2}{\kappa_{\beta_1}-1} \left(-\frac{1}{24}\right) 
    (\L \cdot \beta_1)^2 (\beta_2 \cdot \L)^3 (\beta_1 \cdot \beta_2) 
    N^{(0)}_{\beta_1} N^{(0)}_{\beta_2} \nonumber \\
  & +4 \left(- \frac{1}{24}\right) 
    (\L \cdot \L) (\beta \cdot \L)^3 N^{(0)}_{\beta}.
\end{align}

The remaining cycles all have 2 genus zero components so the calculations are 
simpler. We will first consider $\delta_{0,3}$:
\begin{align*}
  N^{\pi^* \delta_{0,3}}_{\beta,X} = \frac{1}{2}
    \sum_{\substack{\beta_1 + \beta_2 = \beta \\ A, B}} \ 
    \sum_{i,j,k,l} g^{ij} g^{kl}
    & N^{(0)}_{\beta_1,X}(\gamma_i, \gamma_j, \gamma_k, \mu_{\alpha} | \alpha \in A)\\ 
  & \times N^{(0)}_{\beta_2,X}(\gamma_l, \mu_{\alpha} | \alpha \in B), 
\end{align*}
where the first sum is over sets $A, B$ satisfying
\begin{equation*}
  A \sqcup B = \{1,\ldots, \kappa_{\beta}+2\}, \quad |A \cap \{1,2,3,4\}| =1.
\end{equation*}
The factor of $\frac{1}{2}$ appears since the dual graph of a generic curve in $\delta_{0,3}$
has an automorphism of order $2$. Neither $\beta_1 = 0$, nor $\beta_2 = 0$ has any non-trivial 
contribution so it is straight forward to see that
\begin{equation}
  N^{\pi^* \delta_{0,3}}_{\beta,X} = \sum\limits_{ \beta_1 + \beta_2 = \beta} 
    2 \binom{\kappa_{\beta}-2}{\kappa_{\beta_1}-1} 
    (\beta_1\cdot \L) (\beta_2 \cdot \L)^3 
    (\beta_1 \cdot \beta_2) (\beta_1 \cdot \beta_1) 
    N^{(0)}_{\beta_1} N^{(0)}_{\beta_2}.
\end{equation}

The calculation for $\delta_{0,4}$ is a bit more subtle:
\begin{align*}
  N^{\pi^* \delta_{0,4}}_{\beta,X} = \frac{1}{2}
    \sum_{\substack{\beta_1 + \beta_2 = \beta \\ A, B}} \ 
    \sum_{i,j,k,l} g^{ij} g^{kl}
    & N^{(0)}_{\beta_1,X}(\gamma_i, \gamma_j, \gamma_k, \mu_{\alpha} | \alpha \in A)\\ 
  & \times N^{(0)}_{\beta_2,X}(\gamma_l, \mu_{\alpha} | \alpha \in B), 
\end{align*}
where the first sum is over sets $A, B$ satisfying
\begin{equation*}
  A \sqcup B = \{1,\ldots, \kappa_{\beta}+2\}, \quad A \cap \{1,2,3,4\} = \emptyset.
\end{equation*}
Contribution from $\beta_2 = 0$ is $0$. When $\beta_1 = 0$, we must have $A = \emptyset$
which leads to 
\begin{align}
  N^{\pi^* \delta_{0,4}}_{\beta,X} = & \frac{1}{2} 
    \sum_{\beta_1 + \beta_2 = \beta} \binom{\kappa_{\beta}- 2}{\kappa_{\beta_1}-1} 
    (\beta_2 \cdot \L)^4 (\beta_1 \cdot \beta_2) (\beta_1 \cdot \beta_1) 
    N^{(0)}_{\beta_1} N^{(0)}_{\beta_2} \nonumber \\  
  & + \frac{1}{2} (2+ b_2(X)) (\beta \cdot \L)^4 N^{(0)}_{\beta}.
\end{align}

Finally, let us consider the cycle $\delta_{0,0}$:
\begin{align*}
  N^{\pi^* \delta_{0,0}}_{\beta,X} = \frac{1}{2}
    \sum_{\substack{\beta_1 + \beta_2 = \beta \\ A, B}} \ 
    \sum_{i,j,k,l} g^{ij} g^{kl}
    & N^{(0)}_{\beta_1,X}(\gamma_i, \gamma_k, \mu_{\alpha} | \alpha \in A)\\ 
  & \times N^{(0)}_{\beta_2,X}(\gamma_j, \gamma_l, \mu_{\alpha} | \alpha \in B), 
\end{align*}
where the first sum is over sets $A, B$ satisfying
\begin{equation*}
  A \sqcup B = \{1,\ldots, \kappa_{\beta}+2\}, \quad |A \cap \{1,2,3,4\}| = 2.
\end{equation*} 
By an analogous calculation as the previous situations we have
\begin{equation}
  N^{\pi^* \delta_{0,0}}_{\beta,X} = \frac{3}{2} 
    \sum_{ \beta_1 + \beta_2 = \beta} \binom{\kappa_{\beta}-2}{\kappa_{\beta_1}-1} 
    (\beta_1 \cdot \L)^2 (\beta_2 \cdot \L)^2 (\beta_1 \cdot \beta_2)^2 
    N^{(0)}_{\beta_1} N^{(0)}_{\beta_2}.
\end{equation}

Now collecting all these terms and using relation \eqref{eq:relation} we obtain the desired
formula \eqref{rec_formula}.

\section{Low degree checks} 
\label{ld_check}
We will  now describe some  concrete low degree checks that we have performed. 
Let $X_k$ be a del-Pezzo surface obtained by blowing up $\mathbb{P}^2$ at $k \leq 8$ 
points. It is a classical fact that
\[ N^{(1)}_{d L + \sigma_1 E_1 + \ldots +\sigma_r E_r, X_k} \,
=\, N^{(1)}_{dL + \sigma_1 E_1 + \ldots + \sigma_{r-1}E_{r-1}, X_{k-1}}, \] 
if $\sigma_r$ is $-1$ (Qi~\cite{Qi}) or $0$ (Theorem 1.3 of Hu~\cite{Hu}). 
We give  a self contained reason for this assertion in our special case. Consider $X_1$ 
which is $\pp^2$ blown up at the point $p$. 
Let us consider the number $N^{(1)}_{dL -E_1, X_1}$; this is the number of genus one curves 
in $X_1$ representing the class $dL-E_1$ and passing through $3d-1$ generic points. 
Let $\mathcal{C}$ be one of the curves counted by the above number. The curve $\mathcal{C}$
intersects the exceptional divisor exactly at one point. 
Furthermore, since the $3d-1$ points are generic, they can be chosen not to lie in the 
exceptional divisor; let us call the points $p_1, p_2, \ldots, p_{3d-1}$. 
Hence, when we consider the blow down from $X_1$ to $\mathbb{P}^2$, the curve $\mathcal{C}$ 
becomes a curve in $\mathbb{P}^2$ passing through $p_1, p_2, \ldots, p_{3d-1}$ and the blow 
up point $p$. We thus get a genus one, degree $d$ curve in $\mathbb{P}^2$ passing through 
$3d$ points. There is a one to one correspondence between curves representing the class 
$dL-E_1$ in $X_1$ passing through $3d-1$ points and degree $d$ curves in $\mathbb{P}^2$ 
passing through $3d$ points. Hence $N^{(1)}_{dL-E_1, X_1} = N^{(1)}_{dL, \pp^2}$. A 
similar argument holds when there are more than one blowup points. The same argument also 
shows that $N^{(1)}_{dL+0E_1, X_1} = N^{(1)}_{dL, \pp^2}$; the same reasoning holds by taking 
a curve in the blowup and then considering its image under the blow down. The blow down gives 
a one to one correspondence between the two sets and hence, the corresponding numbers 
are 	the same. 

We have verified this assertion in many cases. For instance 
we have verified that 
\[N^{(1)}_{5L-E_1- E_2, X_2} = N_{5L -E_1, X_1,}^{(1)} \, 
=\, N_{5L+0E_1, X_1}^{(1)} \,=\, N_{5L, \pp^2}^{(1)} \, . \]
The reader is invited to use our program and verify these assertions. Hence 
without ambiguity we write $N^{(1)}_{dL + \sigma_1 E_1 + \ldots + \sigma_r E_r}$ for
$N^{(1)}_{dL + \sigma_1 E_1 + \ldots + \sigma_r E_r, X_r}$. 

Next, we note that in \cite{Dub}, Dubrovin has computed the genus one Gromov-Witten 
Invariants of $\mathbb{P}^1 \times \mathbb{P}^1$; our numbers agree with the numbers 
he has listed in his paper (Page 463).  

Finally, in \cite{Vak}, Ravi Vakil has explicitly computed some $N^{(g)}_{\beta}$ 
for del-Pezzo surfaces (Page 78). Our numbers agree with the following numbers he has listed:  
\begin{align*}
N^{(1)}_{5L-2E_1}&= 13775, ~~~~N^{(1)}_{5L-2E_1-2E_2-2E_3} = 225, ~~N^{(1)}_{5L-2E_1- 2E_2- 2E_3- 2E_4} = 20, \\ 
N^{(1)}_{5L-3E_1} & = 240, ~~N^{(1)}_{5L-3E_1-2E_2} = 20 \quad \textnormal{and} \quad N^{(1)}_{5L-3E_1-2E_2-2E_3} = 1.  
\end{align*}

\section*{Acknowledgements}
We would like to thank Ritwik Mukherjee for several fruitful discussions.
The second author is indebted to Ritwik Mukherjee specially for suggesting the project 
and spending countless hours of time for discussions. The first author is also very grateful 
to ICTS for their hospitality and conducive atmosphere for doing mathematics research; he would 
specially like to acknowledge the program Integrable Systems in Mathematics, Condensed Matter 
and Statistical Physics (Code: ICTS/integrability2018/07) where a significant part of the project 
was carried out. We are also grateful  to Ritwik Mukherjee for mentioning our result 
in the program Complex Algebraic Geometry (Code: ICTS/cag2018), which was also organized by ICTS. 
The first author was supported by the DST-INSPIRE grant IFA-16 MA-88 during the course of this research.
Finally, the second author would like to thank Ritwik Mukherjee for supporting this project through 
the External Grant he has obtained, namely MATRICS (File number: MTR/2017/000439) that has been 
sanctioned by the Science and Research Board (SERB).

\end{document}